\documentclass[11pt]{amsart}
\usepackage{amsfonts,amsmath,amsthm,amscd,amssymb,latexsym,cite,verbatim,texdraw,floatflt,
caption2}
\RequirePackage{hyperref}
\usepackage[a4paper]{geometry}

\title{Coarse selectors  of groups    }
\author{ Igor  Protasov}

\address{I.Protasov: Taras Shevchenko National University of Kyiv, Department of Computer Science and Cybernetics, Academic Glushkov pr. 4d, 03680 Kyiv, Ukraine}
\email{i.v.protasov@gmail.com}

\begin{document}
\begin{abstract} 
For a group $G$, $\mathcal{F}_G$ denotes the set of all non-empty finite subsets of $G$. We extend the  finitary coarse structure of $G$ from $G\times G$ to $\mathcal{F}_G\times \mathcal{F}_G$  and say that a macro-uniform mapping $f: \mathcal{F}_G \rightarrow 
\mathcal{F}_G$ (resp. $f: [G]^2 \rightarrow G$) is a finitary selector (resp. 2-selector) of $G$ if  $f(A)\in A$ for each $A\in \mathcal{F}_G$ (resp. $ A \in [G]^2 $). We prove that a  group $G$  admits a finitary selector iff $G$ admits a 2-selector and iff $G$ is a finite extension of an infinite cyclic subgroup or $G$ is countable and locally finite. We use this result to characterize groups admitting linear orders compatible with finitary coarse structures.
\end{abstract}
\maketitle



\section{ Introduction and results}

%\normalsize \cite{b1}.

The notions of selectors went from {\it Topology}.
Let $X$ be a topological space, $exp \  X$ 
denotes  the set of all non-empty closed subsets of $X$ endowed with some (initially, the Vietoris) topology,
$\mathcal{F}$ be a non-empty subset of $exp \  X$.
A continuous mapping  
$f: \mathcal{F}\rightarrow X$ is called an 
$\mathcal{F}$-selector of $X$ if $f(A)\in A$ for each 
$A\in \mathcal{F}$. The question on selectors of  
topological spaces was studied in a plenty of papers, we mention only \cite{b1}, \cite{b4}, \cite{b9}, \cite{b10}.

Formally, coarse spaces, introduced independently and simultaniously in \cite{b17} and \cite{b13}, can be 
considered as asymptotic counterparts of uniform topological spaces. 
But actually, this notion is rooted in {\it Geometry, Geometric Group Theory} and {\it Combinatorics}, 
see [17, Chapter 1], [6, Chapter 4] and 
\cite{b13}.  Every group $G$ admits the natural finitary coarse structure which, in the case of finitely generated $G$, can be viewed as the metric 
structure of a Cayley graph of $G$.
At this point, we need some basic definitions.

\vspace{3 mm}

Given a set $X$, a family $\mathcal{E}$  of subsets of $X\times X$ is called a
{\it  coarse structure} on $X$ if

\begin{itemize}
\item{} each $E \in \mathcal{E}$  contains the diagonal $\bigtriangleup _{X}:=\{(x,x): x\in X\}$ of $X$;
\vspace{3 mm}

\item{}  if  $E$, $E^{\prime} \in \mathcal{E}$  then  $E \circ E^{\prime} \in \mathcal{E}$  and
$ E^{-1} \in \mathcal{E}$,    where  $E \circ E^{\prime} = \{  (x,y): \exists z\;\; ((x,z) \in E,  \ (z, y)\in E^{\prime})\}$,    $ E^{-1} = \{ (y,x):  (x,y) \in E \}$;
\vspace{3 mm}

\item{} if $E \in \mathcal{E}$ and  $\bigtriangleup_{X}\subseteq E^{\prime}\subseteq E$  then  $E^{\prime} \in \mathcal{E}$.
\end{itemize}

Elements $E\in\mathcal E$ of the coarse structure are called {\em entourages} on $X$.

For $x\in X$  and $E\in \mathcal{E}$ the set $E[x]:= \{ y \in X: (x,y)\in\mathcal{E}\}$ is called the {\it ball of radius  $E$  centered at $x$}.
Since $E=\bigcup_{x\in X}( \{x\}\times E[x]) $, the entourage $E$ is uniquely determined by  the family of balls $\{ E[x]: x\in X\}$.
A subfamily ${\mathcal E} ^\prime \subseteq\mathcal E$ is called a {\em base} of the coarse structure $\mathcal E$ if each set $E\in\mathcal E$ is contained in some $E^\prime \in\mathcal E^\prime$.

The pair $(X, \mathcal{E})$  is called a {\it coarse space}  \cite{b17} or  a {\em ballean} \cite{b13}, \cite{b16}. 

%We note that  coarse spaces can be considered as asymptotic counterparts of uniform spaces, see [13, Section 1.1].

A coarse space $(X, \mathcal{E})$  is called
 {\it connected} if,  
 for any $x, y \in X$, there exists $E\in \mathcal{E}$ such that $y\in E[x]$.

A subset  $Y\subseteq  X$  is called {\it bounded} if $Y\subseteq E[x]$ for some $E\in \mathcal{E}$
  and $x\in X$.
If  $(X, \mathcal{E})$  is connected then 
the family $\mathcal{B}_{X}$ of all bounded subsets of $X$  is a bornology on $X$.
We recall that a family $\mathcal{B}$  of subsets of a set $X$ is a {\it bornology}
if $\mathcal{B}$ contains the family $[X] ^{<\omega} $  of all finite subsets of $X$
 and $\mathcal{B}$  is closed   under finite unions and taking subsets. A bornology $\mathcal B$ on a set $X$ is called {\em unbounded} if $X\notin\mathcal B$.
A subfamily  $\mathcal B^{\prime}$ of $\mathcal B$ is called a base for $\mathcal B$ if, for each $B \in \mathcal B$, there exists $B^{\prime} \in \mathcal B^{\prime}$ such that $B\subseteq B^{\prime}$.

Each subset $Y\subseteq X$ defines a {\it subspace}  $(Y, \mathcal{E}|_{Y})$  of $(X, \mathcal{E})$,
 where $\mathcal{E}|_{Y}= \{ E \cap (Y\times Y): E \in \mathcal{E}\}$.
A  subspace $(Y, \mathcal{E}|_{Y})$  is called  {\it large} if there exists $E\in \mathcal{E}$
 such that $X= E[Y]$, where $E[Y]=\bigcup _{y\in Y} E[y]$.

Let $(X, \mathcal{E})$, $(X^{\prime}, \mathcal{E}^{\prime})$
 be  coarse spaces. 
 A mapping $f: X \to X^{\prime}$ is called
  {\it  macro-uniform }  if for every $E\in \mathcal{E}$ there
  exists $E^{\prime}\in \mathcal{E}^{\prime}$  such that $f(E(x))\subseteq  E^{\prime}(f(x))$
    for each $x\in X$.
If $f$ is a bijection such that $f$  and $f ^{-1 }$ are macro-uniform, then   $f  $  is called an {\it asymorphism}.
If  $(X, \mathcal{E})$ and  $(X^{\prime}, \mathcal{E}^{\prime})$  contain large  asymorphic  subspaces, then they are called {\it coarsely equivalent.}

Given a coarse  spaces
$(X, \mathcal{E})$, we denote by 
$exp \ X$ the set of all non-empty subsets of $X$ and endow  $exp \ X$ with the coarse structure $exp \ \mathcal{E}$ with the base 
$\{ exp \ E: E\in \mathcal{E} \}$,
where 
$$(A,B)\in exp \ E 
  \Leftrightarrow A \subseteq E[B], \ \ B\subseteq E[A].$$
The coarse space $(exp \ X, exp \ \mathcal{E} )$ is called the {\it hyperballean} of 
$(X,\mathcal{E})$, 
for hyperballeans see \cite{b2}, \cite{b3}, \cite{b14}, \cite{b15}.
  
Now we are ready to the key definition.
Let   $(X, \mathcal{E})$ be coarse  space, 
  $\mathcal{F}$ be 
  a non-empty  subspace  of $exp \ X$. 
  A macro-uniform mapping 
 $f: \mathcal{F} \longrightarrow X$ 
 is called 
 an $\mathcal{F}$-{\it selector} 
 of $(X,\mathcal{E})$ if $f(A)\in A$ for each $A\in \mathcal{F}$. In the case 
 $\mathcal{F}= exp \ X$,
 $\mathcal{F}= \mathcal{B}\setminus \{0\} $,
 $\mathcal{F}= [X]^2$ we get a {\it global selector},
 a {\it bornologous selector}  
 and a {\it 2-selector}  respectively.
 The investigation of selectors of coarse was  initiated in \cite{b11}, \cite{b12}.

 Every group $G$ with  the identity $e$ can be considered as the 
 coarse  spaces
$(G, \mathcal{E})$, where $\mathcal{E}$ is the 
{\it (right) finitary coarse structure} with the base 
 $$\{\{ (x,y) : x\in Fy \} : F\in [G]^{<\omega}, \  e\in F \} .$$

 We note that the bornology of $(G, \mathcal{E})$ coincides with $\mathcal{F}_G$ and use the name {\it finitary selector} in place the bornologous selector.

 Every metric $d$ on a set $X$ defines the coarse structure $\mathcal{E}_d$ on $X$ with the base 
 $\{\{ (x,y) : d(x,y)\leq r \} : r> 0 \}$.
 Given a connected graph $\Gamma$, $\Gamma =\Gamma [V]
 $, we denote by $d$ the path metric on the set 
 $V$ of vertices of $\Gamma$ and consider  $\Gamma$ as the 
 coarse space $(V, \mathcal{E}_d)$.
 We recall that $\Gamma$ is  {\it locally finite}
 if the set $\{ y: d(x,y)\leq 1\}$ if finite for each 
 $x\in V$.
 
 Our goal is to prove  the following theorem.
 
 \vspace{7 mm}

 {\bf Theorem 1. } 
{\it For a group $G$,
 the following statements are equivalent
 \vspace{5 mm}
 
 $(i)$  $G$ admits a finitary  selector; 
  
 \vspace{5 mm}
 
 $(ii)$  $G$ admits a 2-selector; 
  
\vspace{5 mm}
 
 $(iii)$  $G$ is 
  a finite extension of an
  infinite  cyclic subgroup or 
 $G$ is countable and locally  finite (i.e. every finite subset of $G$ generates a finite subgroup).
 }
  \vspace{5 mm}
  
 In the prof of  Theorem 1
  we use the following characterization of locally finite graphs admitting 
 selectors. 
 By $\mathbb{N}$ and $\mathbb{Z}$,
 we denote  graphs, on the sets of natural and integer 
 numbers in which two vertices  $a, b$ are incident if and only if  $|a-b|=1$. 
 We note also that two graphs are coarsely equivalent  if and only if they are quasi-isometric, see [6, Chapter 4] for quasi-isometric spaces.

\vspace{7 mm}

 {\bf Theorem 2. } 
{\it For a locally finite graph  $\Gamma$,
 the following statements are equivalent:
 \vspace{3 mm}
 
 $(i)$  $\Gamma$ admits a finitary  selector; 
  
 \vspace{3 mm}
 
 $(ii)$  $\Gamma$ admits a 2-selector; 
  
\vspace{3 mm}
 
 $(iii)$  $\Gamma$ is either 
   finite   or    
 coarsely equivalent to   
    $\mathbb{N}$ and $\mathbb{Z}$.}
 
 \vspace{5 mm}
 We  prove   Theorem 2
in Section 2 and  Theorem 1 in Section 3.
In Section 4, we apply Theorem 1 to characterize groups 
admitting linear orders compatible with finitary coarse structures.

\section{ Proof of Theorem 2 }

The implication $(i)\Rightarrow (ii)$ is evident.
To prove $(ii)\Rightarrow (iii)$,  
we choose a 2-selector $f$ of $\Gamma [V]$
and get $(iii)$ at the end of some chain  of elementary observations. 

We define a binary relation $\prec$ on $V$  as follows: $a\prec b$ iff $a\neq b$ and $f(\{a, b\})=a$.

We use also the Hausdorff metric on the set of all
non-empty
 finite subsets of $V$ defined by
 $d_H (A,B)= max \{ d(a,B)$, $ \ d(b, A): a\in A, b\in B\}$,
$d(a,B)= min \{ d(a,b): b\in B\}$.
We note that the coarse structure on 
$[V]^2$ is defined by $d_H$.
Since $f$ is macro-uniform, there exists the minimal natural number $r$ such that if $A,B\in[V]^2$
and $d_H(A, B)\leq 1$ then $d (f(A), f(B))\leq r$.
We fix and use this $r$.

We recall that a sequence of vertices 
$a_0, \dots , a_m$ is a {\it geodesic path} if $d(a_0, a_m)= m$
and $d(a_i, a_{i+1})= 1$ for  each $i\in \{0, \dots , m-1 \}$.
\vspace{7 mm}

{\bf Claim 1.} {\it Let 
$a_0, \dots , a_m$ be a geodesic path in $V$ and $m\geq r$.
If $a_0 \prec  a_r$   (resp. $a_r \prec  a_0$) then 
$a_i \prec  a_j$ (resp. $a_j \prec  a_i$) for all 
$i,j$ such that $j-i \geq r$.}

\vspace{7 mm}

Let $a_0 \prec  a_r$.
By the choice of $r$, we have 
$a_0 \prec  a_{r+1}, \dots a_0 \prec  a_{j}$ and 
$a_1 \prec  a_{j}, \dots a_i \prec  a_{j}$.

\vspace{7 mm}

{\bf Claim 2.} {\it Let $v\in V$,
$B(v, r)= \{ x\in V: d(x, v)\leq  r\}$
and $U$ be a subset of $V\setminus B(v,r)$ such that the graph $\Gamma[U]$ is connected. 
Then either $v\prec u$ for each $u\in U$ or $u\prec v$ for each $u\in U$.}
\vspace{5 mm}

We take arbitrary $u, u^\prime \in U$ and choose $a_0, \dots, a_k$
in $U$ such that $a_0=u$, $a_k=u^\prime$ and 
$d(a_i, a_{i+1})=1$ for each $i\in \{0, \dots, k-1\}$.
Let $a_0\prec v$. By the choice of $r$, we have $a_1\prec v$, $\dots, a_k\prec v$.

\vspace{7 mm}

{\bf Claim 3.} {\it Let $u,v,v^\prime\in V$, 
$d(v, v^\prime)=n$ and  $d(u, v)>n+r$.  If $u\prec v$ (resp. $v\prec u$) then $u\prec v^\prime$ (resp. $v^\prime \prec u$)}.

\vspace{7 mm}

We choose a geodesic path $a_0, \dots , a_m$ from $v$ to $v^\prime$.
Let $u\prec v$. By the choice of $r$,  $u\prec a_0$,
$u\prec a_1, \dots, u\prec a_n$.

\vspace{7 mm}

{\bf Claim 4.} {\it Let $a_0, \dots , a_m$ be a geodesic path in $V$, $v\in V$,
$d(v, \{a_0, \dots , a_m\}) = d(v, a_k)$, $k> 2r+1$,
$m-k> 2r+1$. Then $d(v, a_k)\leq r$.}
\vspace{7 mm}

We take the first alternative given by Claim 1, 
the second is analogical. 
Then 
$a_0\prec a_k, a_k\prec a_m$.
Assuming that $d(v, a_k)> r$, we can replace $v$ to some point on a geodesic path from $v$ to $a_k$ and get 
$d(v, a_k)= r+1$.
We take the first alternative given by Claim 2, the second is analogical. 
Then $v\prec a_0$, $v\prec a_m$.
But  $v\prec a_0$ and  $a_0\prec a_k$
contradict Claim 3.

\vspace{7 mm}

We recall that a sequence $(a_n)_{n<\omega}$  in $V$
is a {\it ray } if $d(a_i , a_j)=j-i$ for all $i<j$. 
Evidently, $\Gamma[\{ a_n : n<\omega\}]$ is  asymorphic to $\mathbb{N}$.

\vspace{7 mm}

{\bf Claim 5.} {\it Let $(a_n)_{n<\omega}$ , $(c_n)_{n<\omega}$  be rays in $V$,
 $A=\{a_n : n<\omega \},$
$C=\{c_n : n<\omega \}$  and  $A\cap C = \emptyset$.
Let $t_0 , \dots , t_k$
be a geodesic path from $a_0$ to $c_0$,
$T=\{ t_0, \dots, t_k\}.$
Assume that $T\cap \{A\} = \{ a_0 \}$, $T\cap C = \{ c_0 \}$.
If there exists a finite subset $H$ of $V$ such that every geodesic path from a vertex $a\in A$ to a vertex  $c\in C$ meets $H$ then $(A\cup C \cup T, d)$ is asymorphic to $\mathbb{Z}$.}

\vspace{7 mm}

We define a bijection 
$f: A\cup C \cup T \rightarrow  \mathbb{Z} $ by
$$f(c_i)= -i-1, \ \  f(t_i)= i, \ \ f(a_i)= i+k + 1 $$
and show that $f$ is an asymorphism. 

If $x,y \in A\cup C \cup T$ then $|f(x)- f(y)|\leq d(x,y)$.
Hence, $f^{-1}$ is macro-uniform.

We denote by $p=max \{ d(a_0 , h), d(b_0 , h): h\in H \}$.
Then the restriction of $f$ to $  C \cup T\cup \{a_0 , \dots , a_p \}$ is an asymorphism and the restriction of $f$ to $ A\cup T \cup \{c_0 , \dots , c_p \}$  
is an asymorphism. Let $n>p$, $m>p$. Since a geodesic path from 
$c_n$ to $a_m$ meets $H$, we have 
$$d(a_m , c_n)\leq n-p+m-p=|f(a_m)- f(c_n)| - k -2p,$$
so $f$ is macro-uniform and the claim is proven.

\vspace{7 mm}
We suppose that $V$ is infinite.
Since 
$\Gamma[V]$ is locally finite,
there exists a ray $(a_n)_{n<\omega}$ in $V$.
We put $A=\{ a_n: n<\omega\}$.
If $V\setminus B(A, r)$ is finite then $\Gamma[V]$ is coarsely equivalent to $\mathbb{N}$.

We suppose $V\setminus B(A, r)$ is infinite, take $u\in V\setminus B(A, r)$ and show that every path $P$ from $u$ to a point from $B(A,r)$ meets 
$B(\{ a_0, \dots , a_{2r+1}\}, r+1)$.
We take a point $v\in P$ such that
$d(v, A)= r+1 $ and take $k$ such that $d(v,a_k)=r+1$.
By Claim 4, $k\leq 2r+1$, so 
$v\in B(\{ a_0, \dots , a_{2r+1}\}, r+1)$.
We choose a ray  $(c_n)_{n<\omega}$ in $V\setminus B(A, r)$ and  put $C=(c_n)_{n<\omega}$.
We delete  (if necessary) a finite number of points from $A$ so that $A, C$ and $T$ satisfy the assumptions of Claim 5 with $F=B(\{ a_0, \dots , a_{2r+1}\}, r+1)$.
Then $(B(A\cup C\cup T), d)$ is coarsely equivalent to 
$\mathbb{Z}$.

We show that $V\setminus B(A\cup C, r)$ is finite, so
$\Gamma[V]$ is coasly equivalent to $\mathbb{Z}$.
We suppose the contrary and choose a ray 
$(x_n)_{n<\omega}$ in $V\setminus B(A\cup C, r)$. 
Applying arguments from above paragraph, we can construct a subset $X$ of $V$ such that $(X, d)$ is coarsely equivalent to a tree $T$ which is a union of three rays with common beginning. Since $(X,d)$ has a 
2-selector, by Proposition 5 from  \cite{b12}, $T$ also
admits a 2-selector.
On the other hand, Claim 4 states that $T$ does not admit a 
2-selector and we get a contradiction.

It remains to prove $(iii)\Rightarrow(i)$.
This is evident if $\Gamma$ is finite.
By [12, Proposition 5], it suffices to show that 
$\mathbb{N}$ and $\mathbb{Z}$ admit finitary selectors.
In both cases, a mapping $f$ defined by 
$f(A)=max \ A$ is finitary
 selector. 

\section{ Proof of Theorem 1 }

 Let $G$ be a group with the finite system  $S$ of generators, $S= S^{-1}$.
We recall that the Cayley graph $Cay(G, S)$ is a graph with the set of vertices $G$ and the set of edges 
$\{(x,y): x\neq y, \  xy^{-1}\in S \}$.
We note that the finitary coarse space of $G$ is asymorphic to the coarse space of $Cay(G, S)$.
\vspace{7 mm}

Now let $G$ be an arbitrary group. The implication $(i)\Rightarrow(ii)$ is evident.

\vspace{7 mm}
We prove $(ii)\Rightarrow(iii)$.
 By 
[12, Theorem 4], $G$ is countable.
Let $f$ be a 
2-selector of $G$.
  We use the binary relation $\prec$ on $G$,
defined in Section 2, and consider two cases.
 
 \vspace{7 mm}
 
 {\it Case 1.} $G$ has an element $a$ of infinite order.
 We denote by $A$ the subgroup of $G$, generated by 
 $a$, and show that $|G: A|$ is finite.
 \vspace{5 mm}
 
 On the contrary, let $|G: A|$ is infinite. We put 
 $S=\{e, a, a^{-1}\}$, denote by $\Gamma[A]$  the graph 
 $Cay(A, S)$ and choose a natural number $r$ such that if $B, C\in [A]^2$ and $d_H (B, C)\leq 1$ then 
 $d (f(A), f(B))\leq r$.
 By Claim 1, either $a^m \prec a^n$ for all $m, n \in \mathbb{Z}$ such that $n-m\geq r$ or $a^n \prec a^m$
 for all $m,n \in \mathbb{Z}$ such that  $n-m \geq r$.
 
 Since $f: [G]^2 \rightarrow G$ is macro-uniform, there exists a finite subset $F$ of $G$ such that $F=F^{-1}$, $e\in F$ and if $B, C\in[G]^2$
 and $A\subseteq SB$, $B\subseteq SA$ then 
 $f(A)\in Ff(B)$.
Since $|G: A|$ is infinite, we can choose $h\in G\setminus FA$, so $Fh\cap A= \emptyset$.
Then either $a^n\prec h$ for each $n\in \mathbb{Z}$ or
$h\prec a^n $ for each $n\in \mathbb{Z}$.
We consider the first alternative, the second is analogical.

Since $f$ is macro-uniform, we can choose $m\in \mathbb{N}$,
$m\geq r$ such that $e\prec a^m $ and  $h\prec a^m $,
but $h\prec a^m $ contradicts above paragraph.

\vspace{7 mm}
 
 {\it Case 2.} $G$ is a torsion group. 
 We suppose that $G$ is not locally finite, choose a finite subset $S$ of $G$ such that the subgroup $H$, 
 generated by $S$, is infinite.
 We denote $\Gamma[H]= Cay (H, S)$.
 By Theorem 2, $\Gamma[H]$ is coarsely equivalent to $\mathbb{N}$ or $\mathbb{Z}$.
 \vspace{5 mm}
 
 We  take $v\in\Gamma[H]$ and denote $S(v, n)= \{ u\in H:  d(v,u) =n \}$, $n\in \mathbb{N}$. By 
  [8, Theorem 	1] or [16, Theorem 5.4.1], there exists a natural number $k$ such that $|S(v,n)|\leq k$ for each $n\in \mathbb{N}$.
  Hence, $H$ is of linear growth. Applying either 
  \cite{b5} or \cite{b7}, we conclude that $H$ has an element of infinite order, a contradiction with the choice of $G$.
  
  \vspace{7 mm}
  
  It remains to verify $(iii)\Rightarrow(i)$.
  If $G$ is a finite extension of an infinite cyclic 
  subgroup then we apply Theorem 2. If $G$ is countable and locally finite, one can refer to Theorem 5 in \cite{b12}, but we  give the following direct proof
  to use in the proof of Theorem 3.
  
  We write $G$ as the union of an increasing chain 
  $\{G_n : n<\omega \}$, $G_0 =\{e\}$ of finite subgroup. 
  For each $n$, we choose some system $R_n$, $e\in R_n$
  of representatives of right cosets of $G_{n+1}$ by
  $G_{n}$, so $G_{n+1}= G_{n} R_n$.
  We denote 
  $$X= \{ (x_n)_{n<\omega} : x_n\in R_n \ and \ x_n =e
  \ for \ all \ but \ finitery \ many \ n \}$$
  and define a bijection $h: G\rightarrow X$ as follows.
  
  We put $h(e)=(x_n)_{n<\omega}, \ x_n =e$.
  Let $g\in G$, $g\neq e$. We choose $n_0$ such that 
  $g\in G_{n_0 +1} \setminus  G_{n_0}$ and write
   $g=g_0 r_{n_0}$, $g_0\in G_{g_0}, \  r_{n_0}\in R_{n_0}$.
   If $g_{0}\neq e$ then we find $n_1$, $g_1 \in G_{n_1}$, $r_{n_1}\in R_{n_1}$ such that $g_0 = g_1 r_{n_1}$.
   After a finite number $k$ of steps, we get $g= r_{n_k} \dots r_{n_1} \ r_{n_0}$.
   We put $h(g)= (y_n)$, where 
   $y_n = r_n$ if $n\in \{ n_k , \dots , n_0 \}$, 
   otherwise, $y_n = e$.
   
   Now we define a linear order $\leq$ on $X$.
   For each $n<\omega$, we choose some linear order 
   $\leq_n$ on $R_n$ with the minimal element $e$.
   If $(x_n)_{n<\omega} \ \neq \ (y_n)_{n<\omega}$
   then we choose the minimal $k$ such that $x_n =y_n$
   for each $n>k$.
   If $x_k <_k  \ y_k $
   then we put  
   $(x_n)_{n<\omega} \ < \ (y_n)_{n<\omega}$.
   
   We note that $(X, \leq)$ is well-ordered, so every non-empty subset of $X$ has the minimal element.
   To define a global selector $f: exp \ G \rightarrow G$, we take an arbitrary $A\in exp \ G$ and put $f(A)= min \ h(A)$.

\section{ Linear orders}

Let $(X, \mathcal{E})$ be a coarse space.
We say that a linear order $\leq$ on $X$ is 
{\it compatible with the coarse structure} 
$\mathcal{E}$
if, for every $E\in \mathcal{E}$, there exists $F\in \mathcal{E}$ such that  
$E\subseteq F$ and if $\{x,y\}\in [X]^2$, $x< y$ 
$(y<x)$
 and $y\in X\setminus F[x]$  then $x^\prime <y$ 
$(y<x^\prime)$ 
  for each $x^\prime \in E[x]$.

\vspace{7 mm}

Let $(X, \mathcal{E})$ be a coarse space, $\leq$ be a linear order on $X$. 
 We say that an entourage $E\in \mathcal{E}$ is interval (with respect to $\leq$) if, for each 
 $x\in X$, there exist $a_x , b_x \in X$ such that 
 $a_x \leq x\leq b_x$ and $E[x]= [a_x, b_x]$.
We say that   
$ \mathcal{E}$  is an interval coarse structure if 
there is a base of $\mathcal{E}$
consisting of interval entourages. Clearly, if 
 $\mathcal{E}$ is interval then $\leq$ is compatible with $\mathcal{E}$.

\vspace{7 mm}

{\bf Theorem 3. } 
{\it Let $G$ be a group,
$\mathcal{E}$ denotes the finitary 
coarse structure on $G$. 
Then the following 
 statements are equivalent
 \vspace{5 mm}
 
 $(i)$  there exists a linear order $\leq$ on $G$
 such that 
 $\mathcal{E}$
 is interval with respect to $\mathcal{E}$; 
  
 \vspace{5 mm}
 
 $(ii)$  there exists a linear order $\leq$ on $G$
 compatible with  
 $\mathcal{E}$;

\vspace{5 mm}
 
 $(iii)$   $G$
 admits a 2-selector. 
 \vspace{7 mm}

Proof.} The imlication $(i)\Rightarrow (ii)$ is  evident, $(ii)\Rightarrow (iii)$ follows from Proposition 2 in 
 \cite{b12}. To prove $(iii)\Rightarrow (i)$, we use 
 Theorem 1 and consider two cases. 
 
 \vspace{7 mm}
 
 {\it  Case 1.} $G$ is a finite extension of an infinite cyclic group $A$.
 We can suppose that $A$ is a normal subgroup.
 Let $A= \{ a^n : n\in \mathbb{Z} \}$, $\{f_0 , \dots , f_m  \}$  be a set of representatives of cosets of  $G$ by $A$, $f_0 = e$, $F=\{f_0 , \dots , f_m  \}$. 
 We set $F_n = F\{ a^{-n}, \dots , a^n \}$,
 $E_n = \{( x,y): xy^{-1}\in F_n\}$ and note that 
 $\{ E_n : n\in \omega \}$ is a base for $\mathcal{E}$.
 
 We endow $G$ with a linear order $\leq$ defined by the  rule: $f_i a^k < f_j a^n $ if and only if either $k< n$
 or $k=n$ and $i< j$.
 
 We choose the minimal natural number $d$ such that 
 $f_j f_j \in F \{ a^{-d}, \dots , a^d \} $ for all 
 $i, j \in \{0, \dots , m \}$.
 Since $f_i ^{-1} a f_j \in  \{a, a^{-1} \}$,
 we have 
 $$E_n [f_i a^k] \ \subseteq \ [f_0 a^{k-n-d}, \ f_m a^{k+n+d}].$$
 
 On the other hand, 
$$[f_0  \ a^{k-n}, \  f_m \ a^{k+n}] \ \subseteq \ 
F_n \ a^{k} \ \subseteq \ F_n \ f_i ^{-1} \ (f_i \ a^k).$$
Hence, $\mathcal{E}$ has an interval base with respect to $\leq$. 
 
  \vspace{7 mm}
 
 {\it Case 2.} $G$ is countable and locally finite. Then $\mathcal{E}$ is interval with respect to the linear order $\leq$ defined in the proof of Theorem 1. 
$ \ \  \  \Box $

\vspace{7 mm}

Let $\leq$ be a linear order on $G$ compatible with $\mathcal{E}$. Does there exist a global selector of $G$? The following theorem gives the negative answer.

\vspace{7 mm}

{\bf Theorem 4. } 
{\it The group $\mathbb{Z}$ does not admit a global selector. 
 \vspace{7 mm}

Proof.} We suppose the contrary and let $f$ be a global selector. Since $f$ is macro-uniform, there exists a natural number $n$ such that if $X, Y \in exp \ G$
and 
$$X\subseteq [-1, 1] + Y, \  Y\subseteq [-1, 1] + X$$
then $f(Y)\in [f(X)-n, \ f(x)+n].$

We put $A=(n+1)\mathbb{Z}$, $a=f(A)$, $A^\prime = A \setminus \{a\}$. 
Then 
\vspace{5 mm}

$f(A^\prime \cup \{a-1\})\in \{ a-1, a- (n+1)\}$, 
$f(A^\prime \cup \{a-2\})\in $ 
$ \{ a-2, a- (n+1)\}, \dots$,  
$f(A^\prime) = a-(n+1)$,
\vspace{5 mm}

$f(A^\prime \cup \{a+1\})\in \{ a+1, a+ n+1\}$, 
$f(A^\prime \cup \{a+2\})\in \{ a+2, a+ n+1\}$, 
$\dots,  f(A^\prime) = a+n+1$,

\noindent and we get a contradiction. 
$ \ \  \  \Box $

\end{document}